\documentclass{article}
\usepackage[T2A]{fontenc}
\usepackage[cp1251]{inputenc}
\usepackage[english,russian]{babel}
\usepackage[tbtags]{amsmath}
\usepackage{amsfonts,amssymb,mathrsfs,amscd}
\numberwithin{equation}{section}
\newtheorem{remark}{Remark}[section]
\newtheorem{lemma}{Lemma}[section]
\newtheorem{theorem}{Theorem}[section]

\newtheorem{condition}{Condition}

\newcommand{\sign}{\mathop{\rm sign}}

\renewcommand{\span}{\mathop{\rm span}}
\renewcommand{\Im}{\mathop{\rm Im}}
\renewcommand{\Re}{\mathop{\rm Re}}

\begin{document}
\begin{Large}
\thispagestyle{empty}
\begin{center}
{\bf Inverse spectral problem for a third-order differential operator with non-local potential\\
\vspace{5mm}
V. A. Zolotarev}\\
\end{center}
\vspace{5mm}

{\small {\bf Abstract.} Spectral problem for a self-adjoint third-order differential operator with non-local potential on a finite interval is studied. Elementary functions that are analogues of sines and cosines for such operators are described. Direct and inverse problems for third-order operators with non-local potential are solved.}
\vspace{5mm}


\vspace{5mm}

\begin{center}
{\bf Introduction}
\end{center}
\vspace{5mm}

Direct and inverse problems for second-order differential operators (Sturm -- Liouville) are well studied and play an important role in mathematical physics and in theory of nonlinear equations \cite{1, 2}. As a rule, these studies are based on transformation operators \cite{1, 2}. Absence of transformation operators for differential operators higher than the second-order does not allow generalizing the formalism of the Sturm -- Liouville operators to this class of operators \cite{3}.

Construction of the Lax L-A pair for nonlinear Camassa-Holm and Dagasperis-Procesi equations leads to a third-order differential operator which is interpreted as an operator describing vibrations of a cubic string \cite{4}. Therefore, the statement of inverse problems is natural for such operators.

This work deals with direct and inverse spectral problems for third-order self-adjoint operators of the form
$$(L_\alpha y)(x)=iy'''(x)+\alpha\int\limits_0^ly(t)\overline{v(t)}dtv(x)$$
where $\alpha\in\mathbb{R}$, $v\in L^2(0,l)$ ($0<l<\infty$), domain of which is formed by functions $y\in W_3^2(0,l)$ such that $y(0)=0$, $y'(0)=y'(l)$, $y(l)=0$. Non-local potentials play an important role in mathematical physics \cite{5}. Works \cite{6, 7, 8} study second-order differential operators with non-local potentials.

The paper consists of 4 sections. First section deals with spectral analysis of the self-adjoint operator $L_0$ of the third derivative. New class of functions (similar to sines and cosines) is introduced and its properties are described. In terms of these functions, eigenfunctions of the operator $L_0$ are described and characteristic function is obtained. Spectrum of the operator $L_0$ is described. Section \ref{s2} studies operator $L_\alpha$ and constructs its eigenfunctions and calculates its characteristic function. Section \ref{s3} describes spectral characteristics of the operator $L_\alpha$. Section \ref{s4} presents solution to the inverse problem and specifies technique of recovery of one-dimensional perturbation by spectral data of the operator $L_\alpha$.

\section{Characteristic function and eigenfunctions of the operator $L_0$}\label{s1}

{\bf 1.1.} Denote by $L_0$ the self-adjoint differential operator
\begin{equation}
(L_0y)(x)\stackrel{\rm def}{=}iD^3y(x)\quad\left(D=\frac d{dx}\right)\label{eq1.1}
\end{equation}
acting in the space $L^2(0,l)$ ($0<l<\infty$). Domain $\mathfrak{D}(L_0)$ of the operator $L_0$ is formed by the functions $y\in W_3^2(0,l)$ that satisfy boundary conditions
\begin{equation}
y(0)=0,\quad y'(0)=y'(l),\quad y(l)=0.\label{eq1.2}
\end{equation}
Let $\{\zeta_k\}_1^3$ be roots of the equation $\zeta^3=1$,
\begin{equation}
\zeta_1=1,\quad\zeta_2=-\frac12+\frac{i\sqrt3}2,\quad\zeta_3=-\frac12-\frac{i\sqrt3}2,\label{eq1.3}
\end{equation}
then
\begin{equation}
\begin{array}{lll}
\zeta_1+\zeta_2+\zeta_3=0,\quad\zeta_1\zeta_2\zeta_3=1,\quad\zeta_1-\zeta_2=i\sqrt3\zeta_3,\quad\zeta_2-\zeta_3=i\sqrt3\zeta_1,\\
\zeta_3-\zeta_1=i\sqrt3\zeta_2.
\end{array}\label{eq1.4}
\end{equation}
The equation
\begin{equation}
iD^3y(x)=\lambda^3y(x)\quad(\lambda\in\mathbb{C})\label{eq1.5}
\end{equation}
has three linearly independent solutions $\{e^{i\lambda\zeta_kx}\}_1^3$ and each solution to equation \eqref{eq1.5} is their linear combination. In particular, solution to the Cauchy problem
\begin{equation}
iD^3y(x)=\lambda^3y(x)+f(x);\quad y(0)=y_0,\quad y'(0)=y_1,\quad y''(0)=y_2\label{eq1.6}
\end{equation}
for $f(x)\equiv0$ is
\begin{equation}
y_0(\lambda,x)=\frac13\sum\limits_k\left(y_0+\frac{y_1}{i\lambda\zeta_k}+\frac{y_2}{(i\lambda\zeta_k)^2}\right)e^{i\lambda\zeta_kx}.\label{eq1.7}
\end{equation}

Hence, by the method of variation of constants, one finds solution to the Cauchy problem \eqref{eq1.6} when $f(x)\not=0$,
\begin{equation}
y(\lambda,x)=\frac13\sum\limits_k\left(y_0+\frac{y_1}{i\lambda\zeta_k}+\frac{y_2}{(i\lambda\zeta_k)^2}\right)e^{i\lambda\zeta_kx}-i\int\limits_0^x\frac13\sum\limits_k\frac{e^{i\lambda\zeta_k(x-t)}}{(i\lambda\zeta_k)^2}f(t)dt.
\label{eq1.8}
\end{equation}
Integral term in \eqref{eq1.8} is a solution to equation \eqref{eq1.6} and satisfies zero initial data $y(0)=0$, $y'(0)=0$, $y''(0)=0$ due to \eqref{eq1.4}.

Instead of the exponents $\{e^{i\lambda\zeta_kx}\}_1^3$, one can take another system of fundamental solutions to equation \eqref{eq1.5}, for example, $c(i\lambda x)$, $s(i\lambda x)$, $d(i\lambda x)$, where
\begin{equation}
c(z)=\frac13\sum\limits_ke^{z\zeta_k};\quad s(z)=\frac13\sum\limits_k\frac1{\zeta_k}e^{z\zeta_k};\quad d(z)=\frac13\sum_k\frac1{\zeta_k^2}e^{z\zeta_k}\label{eq1.9}
\end{equation}
($z\in\mathbb{C}$). Functions \eqref{eq1.9} are similar to cosines and sines for a second-order equation.

\begin{lemma}\label{l1.1}
Entire functions of exponential type \eqref{eq1.9} satisfy the relations:

${\rm(i)}$ $s'(z)=c(z)$, $d'(z)=s(z)$, $c'(z)=d(z)$;

${\rm(ii)}$ $\overline{c(z)}=c(\overline{z})$, $\overline{s(z)}=s(\overline{z})$, $\overline{d(z)}=d(\overline{z})$;

${\rm(iii)}$ $c(z\zeta_2)=c(z)$, $s(z\zeta_2)=\zeta_2s(z)$, $d(z\zeta_2)=\zeta_2^2d(z)$;

${\rm(iv)}$ Euler's formula
$$e^{z\zeta_p}=c(z)+\zeta_ps(z)+\zeta_p^2d(z)\quad(1\leq p\leq3);$$

${\rm(v)}$ the functions \eqref{eq1.9} are solutions to the equation $D^3y=y$ and satisfy the initial data at zero:

\begin{equation}
c(0)=1,\quad c'(0)=0,\quad c''(0)=0;\label{eq1.10}
\end{equation}
$$s(0)=0,\quad s'(0)=1,\quad s''(0)=0;$$
$$d(0)=0,\quad d'(0)=0,\quad d''(0)=1;$$

${\rm(vi)}$ the main identity

$$c^3(z)+s^3(z)+d^3(z)-3c(z)s(z)d(z)=1;$$

${\rm(vii)}$ the summation formulas
$$c(z+w)=c(z)c(w)+s(z)d(w)+d(z)s(w);$$
$$s(z+w)=c(z)s(w)+s(z)c(w)+d(z)d(w);$$
$$d(z+w)=c(z)d(w)+s(z)s(w)+d(z)c(w),$$

${\rm(viii)}$
$$3c^2(z)=c(2z)+2c(-z);$$
$$3s^2(z)=d(2z)+2d(-z),$$
$$3d^2(z)=s(2z)+2s(-z),$$

${\rm(ix)}$
$$s^2(z)-d(z)c(z)=d(-z);$$
$$d^2(z)-s(z)c(z)=s(-z);$$
$$c^2(z)-s(z)d(z)=c(-z);$$

${\rm(x)}$ Taylor's formulas
$$c(z)=1+\frac{z^3}{3!}+\frac{z^6}{6!}+...,$$
$$s(z)=z+\frac{z^4}{4!}+\frac{z^7}{7!}+...,$$
$$d(z)=\frac{z^2}{2!}+\frac{z^5}{5!}+\frac{z^8}{8!}+....$$
\end{lemma}

Proof of the lemma follows from the formulas \eqref{eq1.4}.

The roots $\{\zeta_p\}_1^3$ divide the plane $\mathbb{C}$ into three sectors:
\begin{equation}
S_p=\left\{z\in\mathbb{C}:\frac{2\pi}3(p-1)<\arg z<\frac{2\pi p}3\right\}\quad(1\leq p\leq3).\label{eq1.11}
\end{equation}
Relations (iii) \eqref{eq1.10} imply that it is sufficient to know functions \eqref{eq1.9} in one of the sectors $S_p$.

\begin{lemma}\label{l1.2}
Zeros of the functions $c(z)$, $s(z)$, $d(z)$ \eqref{eq1.9} lie on the rays formed by bisectors of the sectors $S_p$ \eqref{eq1.11} and they correspondingly are given by
$$\left\{-\zeta_2^lx_c(k)\right\}_{k=1}^\infty;\quad\left\{-\zeta_2^lx_s(k)\right\}_{k=1}^\infty;\quad\left\{-\zeta_2^lx_d(k)\right\}_{k=1}^\infty$$
where $l=-1$, $0$, $1$ and $x_c(k)$, $x_s(k)$, $x_d(k)$ are nonnegative numbers in ascending order.

Besides, $x_c(k)$ are simple positive roots of the equation
\begin{equation}
\cos\frac{\sqrt{3}}2x=-\frac12e^{-\frac32x}\quad(x_c(1)>0)\label{eq1.12}
\end{equation}
and the numbers $x_s(k)$ and $x_d(k)$ are nonnegative simple roots of the equations
\begin{equation}
\sin\left(\frac{\sqrt3}2x+\frac\pi6\right)=\frac12e^{-\frac32x}\, (x_s(1)=0);\quad\sin\left(\frac{\sqrt3}2x-\frac\pi6\right)=-\frac12e^{-\frac32x}\, (x_d(1)=0).\label{eq1.12'}
\end{equation}

The sequence $x_d(k)$ is interlacing with the sequence $x_s(k)$, which, in its turn, is interlacing with the sequence $x_c(k)$.
\end{lemma}

Asymptotic behavior of $x_c(k)$, $x_s(k)$, $x_d(k)$ as $k\rightarrow\infty$ can be obtained from the equations \eqref{eq1.12}, \eqref{eq1.12'}.

Formulas \eqref{eq1.7}, \eqref{eq1.8} in terms of the functions \eqref{eq1.9} are given by
\begin{equation}
y_0(\lambda,x)=y_0c(i\lambda x)+y_1\frac{s(i\lambda x)}{i\lambda}+y_2\frac{d(i\lambda x)}{(i\lambda)^2};\label{eq1.7'}
\end{equation}
\begin{equation}
y(\lambda,x)=y_0c(i\lambda x)+y_1\frac{s(i\lambda x)}{i\lambda}+y_2\frac{d(i\lambda x)}{(i\lambda)^2}+i\int\limits_0^x\frac{d(i\lambda(x-t))}{\lambda^2}f(t)dt.\label{eq1.8'}
\end{equation}
\vspace{5mm}

{\bf 1.2.} Find eigenfunctions of the operator $L_0$ \eqref{eq1.1}, \eqref{eq1.2}. Relation \eqref{eq1.7'} implies that the function
\begin{equation}
y_0(\lambda,x)=y_1\frac{s(i\lambda x)}{i\lambda}+y_2\frac{d(i\lambda x)}{(i\lambda)^2}\label{eq1.13}
\end{equation}
is the solution to equation \eqref{eq1.5} and satisfies the first boundary condition $y(0)=0$. Since
$$y'_0(\lambda,x)=y_1c(i\lambda x)+y_2\frac{s(i\lambda x)}{(i\lambda)},$$
due to (i) \eqref{eq1.10}, then the boundary condition $y'(0)=y'(l)$ gives
\begin{equation}
y_1=y_1c(i\lambda l)+y_2\frac{s(i\lambda l)}{i\lambda}.\label{eq1.14}
\end{equation}

\begin{remark}\label{r1.1}
If $\lambda=0$, then \eqref{eq1.14} implies that $y_2=0$, thus, $y_0(0,x)=y_1x$, using the third boundary condition \eqref{eq1.2}, one obtains that $y_1=0$, i. e., $y_0(0,x)=0$ and thus the eigenfunction $y_0(0,x)$ is trivial.
\end{remark}

Assuming that $\lambda\not=0$, from \eqref{eq1.14} one finds that
$$y_1=\frac{y_2}{i\lambda}\cdot\frac{S(i\lambda l)}{1-c(i\lambda l)}$$
and thus function \eqref{eq1.13} is
\begin{equation}
y_0(\lambda,x)=-\frac{y_2}{\lambda^2(1-c(i\lambda l))}\{s(i\lambda x)s(i\lambda l)+d(i\lambda x)(1-c(i\lambda l))\}.\label{eq1.15}
\end{equation}
The third boundary condition $y(l)=0$ implies
\begin{equation}
\frac{s^2(i\lambda l)+d(i\lambda l)(1-c(i\lambda l))}{\lambda^2(1-c(i\lambda l))}=0\quad(\lambda\not=0).\label{eq1.16}
\end{equation}

\begin{remark}\label{r1.2}
Numerator and denominator of fraction \eqref{eq1.16}, as $\lambda\not=0$, do not vanish simultaneously. If this is not the case, then $c(i\lambda l)=1$ and $s(i\lambda l)=0$ and thus $d(i\lambda l)=0$, due to ${\rm (vi)}$ \eqref{eq1.10}. This and ${\rm (iv)}$ \eqref{eq1.10} imply that $e^{i\lambda l\zeta_k}=0$ ($k=1$, $2$, $3$), which is impossible.
\end{remark}

The function
\begin{equation}
\Delta(0,\lambda)\stackrel{\rm def}{=}\{s^2(i\lambda l)+d(i\lambda l)(1-c(i\lambda l))\}\label{eq1.17}
\end{equation}
is said to be the {\bf characteristic function of the operator} $L_0$ \eqref{eq1.1}, \eqref{eq1.2}.

\begin{remark}\label{r1.3}
Equality \eqref{eq1.14} and boundary condition $y(l)=0$ for $y_0(\lambda,x)$ \eqref{eq1.13} lead to the system of linear equations
\begin{equation}
\left\{
\begin{array}{ccc}
{\displaystyle y_1(c(i\lambda l)-1)+y_2\frac{s(i\lambda l)}{i\lambda l}=0;}\\
{\displaystyle y_1\frac{s(i\lambda l)}{i\lambda}-y_2\frac{d(i\lambda l)}{\lambda^2=0}}
\end{array}\right.\label{eq1.18}
\end{equation}
relative to $y_1$, $y_2$. Solution to system \eqref{eq1.18} is non-trivial if its determinant $D(\lambda)=\Delta(0,\lambda)$ \eqref{eq1.17} vanishes, $D(\lambda)=0$.
\end{remark}

Using (ix) \eqref{eq1.10}, one writes characteristic function $\Delta(0,\lambda)$ \eqref{eq1.17} as
\begin{equation}
\Delta(0,\lambda)=\frac1{\lambda^2}(d(i\lambda l)+d(-i\lambda l))=\frac2{3\lambda^2}\sum_k\zeta_k\cos\lambda\zeta_kl.\label{eq1.19}
\end{equation}

\begin{lemma}\label{l1.4}
The function $\Delta(0,\lambda)$ \eqref{eq1.20} is an even entire function of exponential type and
\begin{equation}
\overline{\Delta(0,\lambda)}=\Delta(0,\overline{\lambda});\quad\Delta(0,\zeta_2\lambda)=\Delta(0,\lambda).\label{eq1.20}
\end{equation}
Zeros of $\Delta(0,\lambda)$ are given by $\{\pm\zeta_2^l\lambda_k(0)\}_{k=1}^\infty$ ($l=0$, $1$, $2$) where $\{\lambda_k(0)\}_1^\infty$ are positive simple zeros of $\Delta(0,\lambda)$ on $\mathbb{R}_+$, $0<\lambda_1(0)<\lambda_2(0)<...$, each zero $\lambda_k(0)$ lies in the interval ${\displaystyle\left(\frac\pi l(2k-1),\frac\pi l(2k+1)\right)}$ ($k=1$, $2$, ...) and the asymptotic
\begin{equation}
\lambda_k(0)=\frac{2\pi k}l-\frac\pi{3l}+\delta_k\quad(\delta=o\left(\frac1k\right),k\rightarrow\infty)\label{eq1.21}
\end{equation}
 is true.
\end{lemma}

P r o o f. Due to evenness and symmetry $\Delta(0,\zeta_2\lambda)=\Delta(0,\lambda)$ \eqref{eq1.21}, it is sufficient to find zeros of $\Delta(0,\lambda)$ that are situated on $\mathbb{R}_+$. Since $\Delta(0,0)=-2l^2$, then $\lambda=0$ is not a zero of the function $\Delta(0,\lambda)$ \eqref{eq1.19}. Therefore the equation $\Delta(0,\lambda)=0$ is equivalent to the equality
$$\sum\limits_k\zeta_k\cos\lambda\zeta_kl=0$$
or, taking into account \eqref{eq1.3},
$$\cos\lambda l-\cos\frac{\lambda l}2\ch\frac{\lambda\sqrt{3}l}2-\sqrt3\sin\frac{\lambda l}2\sh\frac{\lambda\sqrt3l}2=0,$$
i. e.,
\begin{equation}\label{eq1.22}
\cos\frac{\lambda l}2\left(\cos\frac{\lambda l}2-\ch\frac{\lambda\sqrt3l}2\right)=\sin\frac{\lambda l}2\left(\sin\frac{\lambda l}2+\sqrt3\sh\frac{\lambda\sqrt3l}2\right).
\end{equation}
It is necessary to find positive zeros $\lambda$ ($>0$) of this equation. Since ${\displaystyle\sin\frac{\lambda l}2+}$ ${\displaystyle\sqrt3\sh\frac{\lambda\sqrt3l}2>0}$, for $\lambda>0$, then equation $\eqref{eq1.22}$ is equivalent to the equality
$$\tan\frac{\lambda l}2=f(\lambda)\quad(\lambda>0),$$
where
$$f(\lambda)=\frac{\displaystyle\cos\frac{\lambda l}2\ch\frac{\lambda\sqrt3l}2}{\displaystyle\sin\frac{\lambda l}2+\sqrt3\sh\frac{\lambda\sqrt3l}2}.$$
$f(\lambda)<0$, as $\lambda>0$, and $f(0)=0$. Using
$$f'(\lambda)=l\frac{\displaystyle1-\sqrt3\sin\frac{\lambda l}2\sh\frac{\lambda\sqrt3l}2-\cos\frac{\lambda l}2\ch\frac{\lambda\sqrt3l}2}{\displaystyle\left(\sin\frac{\lambda l}2+\sqrt3\sh\frac{\lambda\sqrt3l}2\right)^2}$$
and obvious inequalities $\sin x<\sh\sqrt3x$, $\cos x<\ch\sqrt3x$ that are true for all $x>0$, one obtains that
$$f'(\lambda)<l\frac{\displaystyle1-\sqrt3\sin^2\frac{l\lambda}2-\cos^2\frac{\lambda l}2}{\displaystyle\left(\sin\frac{\lambda l}2+\sqrt3\sh\frac{\lambda\sqrt3l}2\right)^2}<0.$$
Consequently, the function $f(\lambda)$ decreases on  $\mathbb{R}_+$ and ${\displaystyle f(\lambda)\rightarrow-\frac1{\sqrt3}}$, as $\lambda\rightarrow\infty$, therefore $f(\lambda)$ monotonously decreases starting from $f(0)=0$ and up to the horizontal asymptote ${\displaystyle y\rightarrow-\frac1{\sqrt3}}$. Hence it follows that in each interval ${\displaystyle\left(\frac\pi l(2k-1),\frac\pi l(2k+1)\right)}$ ($k\in\mathbb{N}$) there lies exactly one root $\lambda_k(0)$ of the equation ${\displaystyle\tan\frac{\lambda l}2=f(\lambda)}$ ($\lambda\in\mathbb{R}_+$) and
\begin{equation}
\lambda_k(0)=\frac{2\pi k}l-\varepsilon_k\quad(0<\varepsilon_k<\frac\pi l).\label{eq1.23}
\end{equation}
Numbers $\lambda_k(0)$ are enumerated in accordance with the intervals ${\displaystyle\left(\frac\pi l(2k-1),\right.}$ ${\displaystyle\left.\frac\pi l(2k+1)\right)}$ ($k\in\mathbb{N}$) and are in ascending order.

Asymptotic \eqref{eq1.21} follows from ${\displaystyle f(\lambda)\rightarrow-\frac1{\sqrt3}}$ ($\lambda\rightarrow\infty$).

Function $\Delta(0,\lambda)$ has no complex zeros lying inside sector $S_1$ \eqref{eq1.11} (as in $S_2$, $S_3$). In other case, if $\Delta(0,\mu)=0$ and $\mu\in S_1$, then $y(\mu,x)$ \eqref{eq1.15} is an eigenfunction of the operator $L_0$ corresponding to the eigenvalue $\mu^3$ which is impossible in view of self-adjointness of $L_0$ since $\mu^3$ is a complex number. $\blacksquare$

\begin{theorem}\label{t1.1}
Spectrum $\sigma(L_0)$ of the operator $L_0$ \eqref{eq1.1}, \eqref{eq1.2} is simple and
\begin{equation}
\sigma(L_0)=\{\pm\lambda_k^3(0):\lambda_k(0)>0,k\in\mathbb{N}\}\label{eq1.24}
\end{equation}
where $\lambda_k(0)$ are positive zeros of characteristic function $\Delta(0,\lambda)$ \eqref{eq1.19} (Lemma \ref{l1.4}). Eigenfunctions corresponding to $\lambda_k(0)$ are
\begin{equation}
u(0,\lambda_k(0),x)=\frac1{u_k}\{s(i\lambda_k(0)x)s(i\lambda_k(0)l)+d(i\lambda_k(0)x)[1-c(i\lambda_k(0)l)]\}\label{eq1.25}
\end{equation}
where $u_k=\|u(0,\lambda_k(0),x)\|_{L^2}$.
\end{theorem}
\vspace{5mm}

{\bf 1.2.} Calculate the resolvent $R_{L_0}(\lambda^3)=(L_0-\lambda^3I)^{-1}$ of the operator $L_0$ \eqref{eq1.1}, \eqref{eq1.2} and let $R_{L_0}(\lambda^3)f=y$, then $L_0y=\lambda^3y+f$ and according to \eqref{eq1.8'} the function
\begin{equation}
y(\lambda,x)=y_1\frac{s(i\lambda x)}{i\lambda}-y_2\frac{d(i\lambda x)}{\lambda^2}+i\int\limits_0^x\frac{d(i\lambda(x-t))}{\lambda^2}f(t)dt\label{eq1.26}
\end{equation}
is a solution to equation \eqref{eq1.6} and satisfies the first boundary condition $y(\lambda,0)=0$. Since
$$y'(\lambda,x)=y_1c(i\lambda x)-iy_2\frac{s(i\lambda x)}\lambda-\int\limits_0^x\frac{s(i\lambda(x-t))}\lambda f(t)dt,$$
then the boundary conditions $y'(0)=y'(l)$ and $y(l)=0$ yield the system of linear equations for $y_1$, $y_2$,
\begin{equation}
\left\{
\begin{array}{ccc}
{\displaystyle y_1(c(i\lambda l)-1)-iy_2\frac{s(i\lambda l)}\lambda=\int\limits_0^l\frac{s(i\lambda(l-t))}\lambda f(t)dt;}\\
{\displaystyle y_1\frac{s(i\lambda l)}{i\lambda}-y_2\frac{d(i\lambda l)}{\lambda^2}=-i\int\limits_0^l\frac{d(i\lambda(l-t))}{\lambda^2}f(t)dt}
\end{array}\right.\label{eq1.27}
\end{equation}
coinciding with system \eqref{eq1.18} as $f=0$. Using the fact that determinant of this system $D(\lambda)$ equals $\Delta(0,\lambda)$ \eqref{eq1.17} (Remark \ref{r1.3}), one obtains that
\begin{equation}
\begin{array}{lll}
{\displaystyle y_1=\frac1{\lambda^3D(\lambda)}\int\limits_0^l\{s(i\lambda l)d(i\lambda(l-t))-d(i\lambda l)s(i\lambda(l-t))\}f(t)dt;}\\
{\displaystyle y_2=\frac i{\lambda^2D(\lambda)}\int\limits_0^l\{[1-c(i\lambda l)]d(i\lambda(l-t))+s(i\lambda l)s(i\lambda(l-t))\}f(t)dt.}
\end{array}\label{eq1.28}
\end{equation}
Hence one finds $y(\lambda,x)$ \eqref{eq1.26},
$$y(\lambda,x)=\frac i{\lambda^4\Delta(0,\lambda)}\left\{s(i\lambda x)\int\limits_0^l[d(i\lambda l)s(i\lambda(l-t))-s(i\lambda l)\right.$$
\begin{equation}
\times d(i\lambda(l-t))]f(t)dt-d(i\lambda x)\int\limits_0^l[(1-c(i\lambda l))d(i\lambda(l-t))+s(i\lambda l)\label{eq1.29}
\end{equation}
$$\left.\times s(i\lambda(l-t))]f(t)dt+\int\limits_0^xd(i\lambda(x-t))[s^2(i\lambda l)+d(i\lambda l)(1-c(i\lambda l))]f(t)dt\right\}.$$
To simplify this expression, one use the statement.

\begin{lemma}\label{l1.5}
For all $z$, the identities

${\rm(a)}$ $d(zl)s(z(l-t))-s(zl)d(z(l-t))=s(-zl)d(-zt)-d(-zl)s(-zt);$

${\rm(b)}$ $(c(zl)-1)d(z(l-t))-s(zl)s(z(l-t))=-(d(zl)+d(-zl))c(-zt)$
\begin{equation}
-s(zl)s(-zt)-(c(zl)-c(-zl))d(-zt);\label{eq1.30}
\end{equation}

${\rm(c)}$ $s(zx)[d(zl)s(z(l-t))-s(zl)d(z(l-t))+d(zx)[(c(zl)-1)d(z(l-t)$
$$-s(zl)s(z(l-t))]=d(z(x-l))d(-zt)-d(zx)d(z(l-t))-d(z(x-t))d(-zl)$$
are true.
\end{lemma}

P r o o f. Using the summation formulas (vii) \eqref{eq1.10}, one obtains
$$d(zl)s(z(l-t))-s(zl)d(z(l-t))=d(zl)[c(zl)s(-zt)+s(zl)c(-zt)$$
$$+d(zl)d(-zt)]-s(zl)[c(zl)d(-zt)+s(zl)s(-zt)+d(zl)c(-zt)]$$
$$=s(-zt)[c(zl)d(zl)-s^2(zl)]+d(-zt)[d^2(zl)-s(zl)c(zl)]$$
$$=-s(-zt)d(-zl)+d(-zt)s(-zl),$$
due to (ix) \eqref{eq1.10}, which proves (a) \eqref{eq1.30}. The second equality (b) \eqref{eq1.30} is similarly proved. Substituting (a), (b) into expression for (c), one obtains
$$s(zx)d(-zt)s(-zl)-s(zx)s(-zt)d(-zl)-d(zx)c(-zt)[d(zl)$$
$$-d(zl)]-d(zx)s(-zt)s(zl)-d(zx)d(-zt)[c(zl)-c(-zl)]$$
$$=-(s(zx)s(-zt)+d(zx)c(-zt))d(-zl)-d(zx)[c(-zt)d(zl)$$
$$+s(-zt)s(zl)+d(-zt)c(zl)]+[s(zx)s(-zl)+d(zx)c(-zl)]d(-tz)$$
$$=-d(z(x-t))d(-zl)-d(zx)d(z(l-t))+d(z(x-l))d(-zt),$$
due to (vii) \eqref{eq1.10}. $\blacksquare$

Substituting (c) \eqref{eq1.30} into \eqref{eq1.29}, one obtains
\begin{equation}
\begin{array}{ccc}
{\displaystyle y(x,\lambda)=\frac i{\lambda^4\Delta(0,\lambda)}\left\{\int\limits_0^l[d(i\lambda(x-l))d(-i\lambda t)-d(i\lambda x)d(i\lambda(l-t))\right.}\\
{\displaystyle \left.-d(i\lambda(x-t))d(i\lambda l)]f(t)dt+\int\limits_0^xd(i\lambda(x-t))[d(i\lambda l)+d(-i\lambda l)]f(t)dt\right\}},
\end{array}\label{eq.1.31}
\end{equation}
hence the statement follows.

\begin{theorem}\label{t1.2}
Resolvent $R_{L_0}(\lambda^3)=(L_0-\lambda^3I)^{-1}$ of the operator $L_0$ \eqref{eq1.2}, \eqref{eq1.3} is
$$(R_{L_0}(\lambda^3)f)(x)=\frac i{\lambda^4\Delta(0,\lambda)}\left\{\int\limits_0^x[d(i\lambda(x-l))d(-i\lambda t)+d(i\lambda(x-t))d(i\lambda l)\right.$$
\begin{equation}
-d(i\lambda x)d(i\lambda(l-t))]+\int\limits_x^l[d(i\lambda(x-l))d(-i\lambda t)-d(i\lambda(x-t))d(-i\lambda l)\label{eq1.32}
\end{equation}
$$-d(i\lambda x)d(i\lambda(l-t))]f(t)dt\}$$
where $\Delta(0,\lambda)$ is the characteristic function \eqref{eq1.20} of the operator $L_0$.
\end{theorem}

\section{Operator $L_\alpha$}\label{s2}

{\bf 2.1} Consider a self-adjoint operator $L_\alpha$ ($=L_\alpha(v)$) in $L^2(0,l)$ which is a one-dimensional perturbation of $L_0$ \eqref{eq1.1}, \eqref{eq1.2},
\begin{equation}
(L_\alpha y)(x)\stackrel{\rm def}{=}iD^3y(x)+\alpha\int\limits_0^ly(t)\overline{v}(t)dtv(x)\label{eq2.1}
\end{equation}
where $\alpha\in\mathbb{R}$ and $v\in L^2(0,l)$. Without loss in generality, assume that $\|v\|_{L^2}=1$. Domains of the operators $L_\alpha$ and $L_0$ coincide, $\mathfrak{D}(L_\alpha)=\mathfrak{D}(L_0)$.

Using the functions $c(z)$, $s(z)$, $d(z)$ \eqref{eq1.9}, define Fourier transforms
$$\widetilde{v}_c(\lambda)\stackrel{\rm def}{=}\langle v(x),c(i\overline{\lambda}x)\rangle=\frac13\sum\limits_k\widetilde{v}_k(\lambda);$$
\begin{equation}
\widetilde{v}_s(\lambda)\stackrel{\rm def}{=}\langle v(x),s(i\overline{\lambda}x)\rangle=\frac13\sum\limits_k\zeta_k^{-1}\widetilde{v}_k(\lambda);\label{eq2.2}
\end{equation}
$$\widetilde{v}_d(\lambda)\stackrel{\rm def}{=}\langle v(x),d(i\overline{\lambda}x)\rangle=\frac13\sum\limits_k\zeta_k\overline{v}_k(\lambda)$$
where
\begin{equation}
\widetilde{v}_k(\lambda)=\int\limits_0^le^{-i\lambda\zeta_kx}v(x)dx\quad(1\leq k\leq3).\label{eq2.3}
\end{equation}
Introduce the operation of involution
\begin{equation}
f^*(\lambda)=\overline{f(\overline{\lambda})},\label{eq2.4}
\end{equation}
then
$$\widetilde{v}_c^*(\lambda)=\langle c(i\lambda x),v(x)\rangle=\frac13\sum\widetilde{v}_k^*(\lambda);$$
\begin{equation}
\widetilde{v}_s^*(\lambda)=\langle s(i\lambda x),v(x)\rangle=\frac13\sum\zeta_k\widetilde{v}_k^*(\lambda);\label{eq2.5}
\end{equation}
$$\widetilde{v}_d^*(\lambda)=\langle d(i\lambda x),v(x)\rangle=\frac13\sum\zeta_k^{-1}\widetilde{v}_k^*(\lambda),$$
besides,
\begin{equation}
\widetilde{v}_1^*(\lambda)=\int\limits_0^le^{i\lambda\zeta_1x}\overline{v}(x)dx;\quad\widetilde v_2^*(\lambda)=\int\limits_0^le^{i\lambda\zeta_3x}\overline{v}(x)dx;\quad\overline v_3^*(\lambda)=\int\limits_0^le^{i\lambda\zeta_2x}\overline{v}(x)dx.\label{eq2.6}
\end{equation}
According to the Cauchy problem \eqref{eq1.6}, equation $L_\alpha y=\lambda^3y$ has the solution (see \eqref{eq1.8'})
\begin{equation}
y(\lambda,x)=y_1\frac{s(i\lambda x)}{i\lambda}-y_2\frac{d(i\lambda x)}{\lambda^2}-i\alpha\langle y,v\rangle\int\limits_0^x\frac{d(i\lambda(x-t))}{\lambda^2}v(t)dt\label{eq2.7}
\end{equation}
satisfying the boundary condition $y(0)=0$. Multiplying equality \eqref{eq2.7} by $\overline{v}(x)$ and integrating it from 0 to $l$, one obtains
\begin{equation}
iy_1\frac{\widetilde{v}_s^*(\lambda)}\lambda+y_2\frac{\widetilde v_d^*(\lambda)}{\lambda^2}+\langle y,v\rangle\left(1+\frac{i\alpha}{\lambda^2}m(x)\right)=0\label{eq2.8}
\end{equation}
where $\widetilde{v}_s^*(\lambda)$ and $\widetilde{v}_d^*(\lambda)$ are given by \eqref{eq1.5} and
\begin{equation}
\begin{array}{lll}
{\displaystyle m(\lambda)\stackrel{\rm def}{=}\left\langle\int\limits_0^xd(i\lambda(x-t))v(t)dt,v(x)\right\rangle=\frac13\sum\zeta_k\psi_k(\lambda);}\\
{\displaystyle\psi_k(\lambda)\stackrel{\rm def}{=}\int\limits_0^ldxe^{i\lambda\zeta_kx}\overline{v}(x)\int\limits_0^xdte^{-i\lambda\zeta_kt}v(t)\quad(1\leq k\leq3).}
\end{array}\label{eq2.9}
\end{equation}

\begin{remark}\label{r2.1}
The functions $\widetilde{v}_k(\lambda)$ \eqref{eq2.3} and $\widetilde{v}_k^*(\lambda)$ \eqref{eq2.6} (and thus the functions \eqref{eq2.2} and \eqref{eq2.5} also) are entire functions of exponential type.
\end{remark}

\begin{lemma}\label{l2.1}
The function $m(\lambda)$ \eqref{eq2.9} is an entire function of exponential type and
\begin{equation}
m(\lambda)+m^*(\lambda)=\widetilde{v}_d(\lambda)\widetilde{v}_c^*(\lambda)+\widetilde{v}_s(\lambda)\widetilde{v}_s^*(\lambda)+\widetilde{v}_c(\lambda)\widetilde{v}_d^*(\lambda).\label{eq2.10}
\end{equation}
\end{lemma}

P r o o f. Since $\psi_k(\lambda)$ \eqref{eq2.9} is the Fourier transform of the convolution
$$\psi_k(\lambda)=\int\limits_0^ldx\int\limits_0^xdte^{i\lambda\zeta_k(x-t)}v(t)\overline{v(x)}=\int\limits_0^ldse^{i\lambda\zeta_ks}\int\limits_s^lv(x-s)\overline{v}(x)dx,$$
then $\psi_k(\lambda)$ is an entire function of exponential type, thus $m(\lambda)$ also has this property.

Relation \eqref{eq2.9} implies that
$$m(\lambda)+m^*(\lambda)=\frac13\{\zeta_1(\psi_1(\lambda)+\psi_1^*(\lambda))+\zeta_2(\psi_2(\lambda)+\psi_3^*(\lambda))+\zeta_3(\psi_3(\lambda)+\psi_2^*(\lambda))\}.$$
Integrating by parts, one obtains
$$\psi_k(\lambda)=\int\limits_0^ldte^{-i\lambda\zeta_kt}v(t)\int\limits_0^ldxe^{i\lambda\zeta_kx}\overline{v(x)}-\int\limits_0^ldte^{-i\lambda\zeta_kt}v(t)\int\limits_0^tdxe^{i\lambda\zeta_kx}\overline{v(x)},$$
therefore
$$\psi_1(\lambda)=\widetilde{v}_1(\lambda)\widetilde{v}_1^*(\lambda)-\psi_1^*(\lambda),\quad\psi_2(\lambda)=\widetilde{v}_2(\lambda)\widetilde{v}_3^*(\lambda)-\psi_3^*(\lambda),$$
$$\psi_3(\lambda)=\widetilde{v}_3(\lambda)\widetilde{v}_2^*(\lambda)-\psi_2^*(\lambda).$$
So,
$$m(\lambda)+m^*(\lambda)=\frac13\{\zeta_1\widetilde{v}_1(\lambda)\widetilde{v}_1^*(\lambda)+\zeta_2\widetilde v_2(\lambda)\widetilde{v}_3^*(\lambda)+\zeta_3\widetilde v_3(\lambda)\widetilde v_2^*(\lambda)\}$$
$$=\int\limits_0^ldx\int\limits_0^ldtd(i\lambda(x-t))\overline{v}(x)v(t)=\int\limits_0^ldx\int\limits_0^ldt\overline{v(x)}v(t)\{c(i\lambda x)d(-i\lambda t)$$
$$+s(i\lambda x)s(-i\lambda t)+d(i\lambda x)c(-i\lambda t)\}$$
due to (vii) \eqref{eq1.10} and according to \eqref{eq2.2}, \eqref{eq2.5},
$$m(\lambda)+m^*(\lambda)=\widetilde{v}_d(\lambda)\widetilde{v}_c^*(\lambda)+\widetilde{v}_s(\lambda)\widetilde{v}_s^*(\lambda)+\widetilde{v}_c(\lambda)\widetilde{v}_d^*(\lambda).\quad\blacksquare$$

Equality \eqref{eq2.10} is an analogue of a well-known statement on the Fourier transform of a convolution.

Formula \eqref{eq2.8} and boundary conditions $y'(0)=y'(l)$, $y(l)=0$ for $y(\lambda,x)$ \eqref{eq2.7} yield the system of linear equations for $y_1$, $y_2$, $\langle y,v\rangle$,
\begin{equation}
\left\{
\begin{array}{lll}
{\displaystyle y_1i\frac{\widetilde{v}_s^*(\lambda)}\lambda+y_2\frac{\widetilde{v}_d^*(\lambda)}{\lambda^2}+\langle y,v\rangle\left(1+\frac{i\alpha}{\lambda^2}m(\lambda)\right)=0;}\\
{\displaystyle y_1(1-c(i\lambda l))+y_2i\frac{s(i\lambda l)}\lambda-\alpha\langle y,v\rangle\int\limits_0^l\frac{s(i\lambda(l-t))}\lambda v(t)dt=0,}\\
{\displaystyle y_1\frac{is(i\lambda l)}\lambda+y_2\frac{d(i\lambda l)}{\lambda^2}+\alpha i\langle y,v\rangle\int\limits_0^l\frac{d(i\lambda(l-t))}{\lambda^2}v(t)dt=0.}
\end{array}\right.\label{eq2.11}
\end{equation}
Solutions to system \eqref{eq2.11} are non-trivial then and only then when its determinant $D_\alpha(\lambda)$ vanishes. The function $\Delta(\alpha,\lambda)=D_\alpha(\lambda)$ is said to be {\bf characteristic function of the operator} $L_\alpha$ \eqref{eq2.1} and equals
\begin{equation}
\begin{array}{ccc}
{\displaystyle\Delta(\alpha,\lambda)\stackrel{\rm def}{=}\left(1+\frac{\alpha i}{\lambda^2}m(\lambda)\right)\Delta(0,\lambda)+\frac{i\alpha}{\lambda^4}\left\{\widetilde{v}_d^*(\lambda)\int\limits_0^l[(c(i\lambda l)-1)d(i\lambda(l-t))\right.}\\
{\displaystyle-s(i\lambda l)s(i\lambda(l-t))]v(t)dt+\widetilde{v}_s^*(\lambda)\int\limits_0^l[d(i\lambda l)s(i\lambda(l-t))}\\
-s(i\lambda l)d(i\lambda(l-t))]v(t)dt\}
\end{array}\label{eq2.12}
\end{equation}
Using \eqref{eq1.19} and (a), (b) \eqref{eq1.30}, one has
$$\Delta(\alpha,\lambda)-\Delta(0,\lambda)=\frac{i\alpha}{\lambda^4}\{m(\lambda)(d(i\lambda l)+d(-i\lambda l))-\widetilde{v}_d^*(\lambda)\widetilde{v}_c(\lambda)$$
$$\times(d(i\lambda l)+d(-i\lambda l))-\widetilde{v}_d^*(\lambda)\widetilde{v}_s(-\lambda)s(i\lambda l)-\widetilde{v}_d^*(\lambda)\widetilde{v}_d(\lambda)(c(i\lambda l)-c(-i\lambda l))$$
$$+\widetilde{v}_s^*(\lambda)\widetilde{v}_d(\lambda)s(-i\lambda l)-\widetilde{v}_s^*(\lambda)v_s(\lambda)d(-i\lambda l)\}.$$

Expressing $\widetilde{v}_s^*(\lambda)\widetilde{v}_s(\lambda)$ from \eqref{eq2.10} and substituting it into this expression, one obtains
$$\Delta(\alpha,\lambda)-\Delta(0,\lambda)=\frac{i\alpha}{\lambda^4}\{m(\lambda)d(i\lambda l)+\widetilde{v}_d(\lambda)[\widetilde{v}_c^*(\lambda)d(-i\lambda l)$$
$$+\widetilde{v}_s^*(\lambda)s(-i\lambda l)+\widetilde{v}_d^*(\lambda)c(-i\lambda l)]-m^*(\lambda)d(-i\lambda l)-\widetilde{v}_d^*(\lambda)$$
$$\times[\widetilde{v}_c(\lambda)d(i\lambda l)+\widetilde{v}_s(\lambda)s(i\lambda l)+\widetilde{v}_d(\lambda)c(i\lambda l)]\}.$$
Taking (vii) \eqref{eq1.10} into account, one finds
$$\widetilde{v}_c(\lambda)d(i\lambda l)+\widetilde{v}_s(\lambda)s(i\lambda l)+\widetilde{v}_d(\lambda)c(i\lambda l)=\int\limits_0^ld(i\lambda(l-x))v(x)dx$$
$$=\int\limits_0^ld(i\lambda t)v(l-t)dt=\widetilde{w}_d(-\lambda)$$
where $\widetilde{w}_d(\lambda)$ is the Fourier transform $\widetilde{w}_d(\lambda)=\langle w(x),d(i\lambda x)\rangle$ \eqref{eq2.2} of the function
\begin{equation}
w(x)\stackrel{\rm def}{=}v(l-x).\label{eq2.13}
\end{equation}

\begin{lemma}\label{l2.2}
Characteristic function $\Delta(\alpha,\lambda)$ \eqref{eq2.12} of the operator $L_\alpha$ \eqref{eq2.1} is
\begin{equation}
\Delta(\alpha,\lambda)-\Delta(0,\lambda)=\frac{i\alpha}{\lambda^4}\{F(\lambda)-F^*(\lambda)\}\label{eq2.14}
\end{equation}
where
\begin{equation}
F(\lambda)\stackrel{\rm def}{=}m(\lambda)d(i\lambda l)+\widetilde{v}_d(\lambda)\widetilde{w}_d^*(-\lambda),\label{eq2.15}
\end{equation}
besides, $\widetilde{v}_d(\lambda)$ and $m(\lambda)$ are given by \eqref{eq2.5}, \eqref{eq2.9}, $\widetilde{w}_d(\lambda)$ is the Fourier transform \eqref{eq2.2} of the function $w(x)$ \eqref{eq2.13}. Moreover, $\Delta(\alpha,\lambda)$ \eqref{eq2.14} is a real entire function of exponential type and
\begin{equation}
\Delta(\alpha,\lambda\zeta_2)=\Delta(\alpha,\lambda),\quad\Delta^*(\alpha,\lambda)=\Delta(\alpha,\lambda).\label{eq2.16}
\end{equation}
\end{lemma}
\vspace{5mm}

{\bf 2.2.} Find eigenfunctions of the operator $L_\alpha$ \eqref{eq2.1}. Boundary conditions $y'(0)=y'(l)$ and $y(l)=0$ for the function $y(\lambda,x)$ \eqref{eq2.7} yield system of equations formed by the two last equalities in \eqref{eq2.11}
$$\left\{
\begin{array}{ccc}
{\displaystyle y_2i\frac{s(i\lambda l)}\lambda-\frac{\alpha\langle y,v\rangle}\lambda\int\limits_0^ls(i\lambda(l-t))v(t)dt=-y_1(1-c(i\lambda l)),}\\
{\displaystyle y_2\frac{d(i\lambda l)}\lambda+\frac{\alpha i\langle y,v\rangle}\lambda\int\limits_0^ld(i\lambda(l-t))v(t)dt=-y_1is(i\lambda l).}
\end{array}\right.$$
Hence we find that
\begin{equation}
\begin{array}{lll}
{\displaystyle\frac{y_2}\lambda=\frac{iy_1}{a(\lambda)}\int\limits_0^l[(c(i\lambda l)-1)d(i\lambda(l-t))-s(i\lambda l)s(i\lambda(l-t))]v(t)dt;}\\
{\displaystyle\frac{\alpha\langle y,v\rangle}\lambda=\frac{y_1}{a(\lambda)}(s(i\lambda l)^2+d(i\lambda l))(1-c(i\lambda l)),}
\end{array}\label{eq2.17}
\end{equation}
besides,
$$a(\lambda)=\int\limits_0^l[d(i\lambda l)s(i\lambda(l-t))-s(i\lambda t)d(i\lambda(l-t))]v(t)dt.$$
Substituting expressions \eqref{eq2.17} into \eqref{eq2.7}, one obtains
$$y(\lambda,x)=-\frac{iy_1}{\lambda a(\lambda)}\left\{\int\limits_0^l\{s(i\lambda x)[d(i\lambda l)s(i\lambda(l-t))-s(i\lambda l)d(i\lambda(l-t))]\right.$$
$$+d(i\lambda x)[(c(i\lambda l)-1)d(i\lambda(l-t)-s(i\lambda l)s(i\lambda(l-t))]\}v(t)dt$$
and according to (c) \eqref{eq1.30},
$$y(\lambda,x)=-\frac{iy_1}{\lambda a(\lambda)}\left\{\int\limits_0^l[d(i\lambda(x-l))d(-i\lambda t)-d(i\lambda x)d(i\lambda(l-t))\right.$$
$$\left.-d(i\lambda(x-t))d(-i\lambda l)]v(t)dt+(d(i\lambda l)+d(-i\lambda l))\int\limits_0^xd(i\lambda(x-t))v(t)dt\right\}$$

\begin{lemma}\label{l2.3}
Eigenfunctions $u(\alpha,\lambda,x)$ of the operator $L_\alpha$ \eqref{eq2.1} are given by
$$u(\alpha,\lambda,x)=\frac1{u(\alpha,\lambda)}\{d(i\lambda(x-l))\widetilde{v}_d(\lambda)-d(i\lambda x)\widetilde{w}_d(-\lambda)$$
\begin{equation}
\left.+d(i\lambda l)\int\limits_0^xd(i\lambda(x-t))v(t)dt-d(-i\lambda l)\int\limits_x^ld(i\lambda(x-t))v(t)dt\right\}\label{eq2.18}
\end{equation}
where $\lambda=\lambda_k(\alpha)$ is a zero of the characteristic function $\Delta(\alpha,\lambda)$ \eqref{eq2.14}, $\widetilde{w}_d(\lambda)$ is the Fourier transform \eqref{eq2.2} of the function $w(x)$ \eqref{eq2.13}; $u(\alpha,\lambda)=\|u(\alpha,\lambda,x)\|_{L^2}$ and
\begin{equation}
u(\alpha,\lambda\zeta_2,x)=u(\alpha,\lambda,x).\label{eq2.19}
\end{equation}
\end{lemma}

\section{Abstract problem}\label{s3}

{\bf 3.1.} Consider an abstract problem of one-dimensional perturbation of a self-adjoint operator with simple spectrum. Let $L_0$ be a self-adjoint operator acting in a Hilbert space $H$ with simple discrete spectrum
\begin{equation}
\sigma(L_0)=\{z_n:z_n\in\mathbb{R},n\in\mathbb{Z}\}.\label{eq3.1}
\end{equation}
Its eigenfunctions $u_n$ ($L_0u_n=z_nu_n$) are orthonormal, $\langle u_n,u_m\rangle=\delta_{n.m}$, and form an orthonormal basis in the space $H$. Resolvent $R_{L_0}(z)=(L_0-zI)^{-1}$ of the operator $L_0$ is
\begin{equation}
R_{L_0}(z)=\sum\limits_h\frac{f_n}{z_n-z}u_n\label{eq3.2}
\end{equation}
where $f_n=\langle f,u_n\rangle$ are Fourier coefficients of the vector $f$ in the basis $\{u_n\}$.

By $L_\alpha$ denote a self-adjoint operator which is a one-dimensional perturbation of $L_0$ (cf. \eqref{eq2.1}),
\begin{equation}
L_\alpha\stackrel{\rm def}{=}L_0+\alpha\langle.,v\rangle v\label{eq3.3}
\end{equation}
where $\alpha\in\mathbb{R}$ and $v$ is a fixed vector from $H$ such that $\|v\|=1$. Domains of $L_\alpha$ and $L_0$ coincide, $\mathfrak{D}(L_\alpha)=\mathfrak{D}(L_0)$.

\begin{lemma}\label{l3.1}
Resolvent $R_{L_\alpha}(z)=(L_\alpha-zI)^{-1}$ of the operator $L_\alpha$ \eqref{eq3.3} is
\begin{equation}
R_{L_\alpha}(z)f=R_{L_0}-\alpha\frac{\langle R_{L_0}(z)f,v\rangle}{1+\alpha\langle R_{L_0}(z)v,v\rangle}\cdot R_{L_0}(z)v\label{eq3.4}
\end{equation}
where $R_{L_0}(z)=(L_0-zI)^{-1}$ is the resolvent of operator $L_0$ and $f\in H$.
\end{lemma}

P r o o f. Let $y=R_{L_0}(z)f$, then
$$f=(L_0-zI)y+\alpha\langle y,v\rangle v,$$
and thus
\begin{equation}
R_{L_0}(z)f=y+\alpha\langle y,v\rangle R_{L_0}(z)v.\label{eq3.5}
\end{equation}
Scalar multiplying this equality  by $v$, one obtains
$$\langle R_{L_0}(z)f,v\rangle=\langle y,v\rangle(1+\alpha\langle R_{L_0}(z)v,v\rangle).$$
From here expressing $\langle y,v\rangle$ and substituting it into \eqref{eq3.5}, one arrives at formula \eqref{eq3.4}. $\blacksquare$

Formulas \eqref{eq3.2}, \eqref{eq3.4} imply
\begin{equation}
R_{L_\alpha}(z)f=\sum\limits_n\frac{f_n}{z_n-z}u_n-\alpha\frac{\displaystyle\sum\limits_k\frac{f_k\overline{v}_k}{z_k-z}}{\displaystyle 1+\alpha\sum\limits_k\frac{|v_k|^2}{z_k-z}}\cdot\sum\limits_n\frac{v_n}{z_n-z}u_n\label{eq3.6}
\end{equation}
where $f_n=\langle f,u_n\rangle$ and $v_n=\langle v,u_n\rangle$ are Fourier coefficients of the vectors $f$ and $v$ correspondingly.

\begin{lemma}\label{l3.2}
If $v_p\not=0$, then point $z_p$ ($\in\sigma(L_0)$ \eqref{eq3.1}) does not belong to the spectrum $\sigma(L_\alpha)$ operator $L_\alpha$ \eqref{eq3.3}.
\end{lemma}

P r o o f. Really, residue of the resolvent $R_{L_\alpha}(z)$ \eqref{eq3.6} vanishes at the point $z_p$,
$$\lim\limits_{z\rightarrow z_p}(z_p-z)R_{L_\alpha}(z)=f_pu_p-\frac{\alpha f_p\overline{v_p}}{\alpha v_p\overline{v_p}}\cdot v_pu_p=0.\quad\blacksquare$$

Hence it follows that it is natural to divide the set $\sigma(L_0)$ \eqref{eq3.1} into two disjoint subsets $\sigma(L_0)=\sigma_0\cup\sigma_1$ where
\begin{equation}
\sigma_0\stackrel{\rm def}{=}\{z_n^0=z_n\in\sigma(L_0):v_n=0\},\quad\sigma_1\stackrel{\rm def}{=}\{z_n^1=z_n\in\sigma(L_0):v_n\not=0\}.\label{eq3.7}
\end{equation}
Assume that elements from $\sigma(L_0)$, $\sigma_0$, $\sigma_1$ are numbered in ascending order. Partition $\sigma(L_0)=\sigma_0\cup\sigma_1$ corresponds to the decomposition of the space $H$ into orthogonal sum $H=H_0\oplus H_1$ where
\begin{equation}
H_0\stackrel{\rm def}{=}\span\{u_n:z_n\in\sigma_0\},\quad H_1\stackrel{\rm def}{=}\span\{u_n:z_n\in\sigma_1\}.\label{eq3.8}
\end{equation}
The subspaces $H_0$, $H_1$ reduce the operator $L_\alpha$ \eqref{eq3.3} and
\begin{equation}
\left.L_\alpha\right|_{H_0}=L_{\alpha,0};\quad\left.L_\alpha\right|_{H_1}=L_{\alpha,1},\label{eq3.9}
\end{equation}
therefore
$$R_{L_\alpha}(z)=R_{L_{\alpha,0}}(z)\oplus R_{L_{\alpha,1}}(z)$$
where
\begin{equation}
\begin{array}{lll}
{\displaystyle R_{L_{\alpha,0}}(z)f=\sum\limits_{z_n^0\in\sigma_0}\frac{f_n}{z_n^0-z}u_n;}\\
{\displaystyle R_{L_{\alpha,1}}(z)f=\frac1{Q(z)}\sum\limits_{z_n^1\in\sigma_1}\frac1{z_n^1-z}\{f_nQ(z)-\alpha v_nP(z,f)\}u_n,}
\end{array}\label{eq3.10}
\end{equation}
besides,
\begin{equation}
Q(z)\stackrel{\rm def}{=}1+\alpha\sum\limits_{z_n^1\in\sigma_1}\frac{|v_n|^2}{z_n^1-z};\quad P(z,f)\stackrel{\rm def}{=}\sum\limits_{z_n^1\in\sigma_1}\frac{f_n\overline{v_n}}{z_n^1-z}.\label{eq3.11}
\end{equation}
Lemma \ref{l3.2} implies that the resolvent $R_{L_{\alpha,1}}(z)$ \eqref{eq3.10} does not have singularities at the points $z_n^1\in\sigma_1$ but can have singularities at zeros of the function $Q(z)$ \eqref{eq3.11}.

The function
\begin{equation}
G(z)=\sum\limits_{z_n^1\in\sigma_1}\frac{|v_n|^2}{z_n^1-z}\label{eq3.12}
\end{equation}
monotonously increases for all $z\in\mathbb{R}\backslash \sigma_1$ (since $F'(z)>0$ for all $z\in\mathbb{R}\backslash \sigma_1$). Therefore equation $1+\alpha F(z)=0$, which is equivalent to $Q(z)=0$, has only simple roots.

\begin{lemma}\label{l3.3}
Zeros of $Q(z)$ \eqref{eq3.11} are real, simple and alternate with numbers $z_n^1\in\sigma_1$ \eqref{eq3.7}.
\end{lemma}

Define the set
\begin{equation}
Q(L_{\alpha,1})\stackrel{\rm def}{=}\{\mu_n:Q(\mu_n)=0\}\label{eq3.13}
\end{equation}
where $Q(z)$ is given by \eqref{eq3.11} and numbers $\mu_n=\mu(\alpha)$ depend on the parameter $\alpha$. Resolvent $R_{L_{\alpha,1}}(z)$ \eqref{eq3.10} cannot have removable singularities at the points $\mu_n\in\sigma(L_{\alpha,1})$. Thus singularities of $R_{L_{\alpha,1}}(z)$ are the simple poles at the points $\mu_n\in\sigma(L_{\alpha,1})$.
\vspace{5mm}

{\bf 3.2} Find eigenfunctions of the operator $L_{\alpha,1}$. Formula \eqref{eq3.10} implies that residue of the resolvent $R_{L_{\alpha,1}}(z)$ at the point $z=\mu_p$ is
\begin{equation}
c_p(f)\stackrel{\rm def}{=}\lim\limits_{z\rightarrow\mu_p}(\mu_p-z)R_{L_{\alpha,1}}(z)=\frac\alpha{Q'(\mu_p)}\sum\limits_{z_n^1\in\sigma_1}\frac{v_n}{z_n^1-\mu_p}\left(\sum\limits_{z_k^1\in\sigma_1}\frac{f_k\overline{v_k}}{z_k^1-\mu_p}\right)u_n.
\label{eq3.14}
\end{equation}
Note that
$$Q'(\mu_p)=\alpha\sum\limits_{z_n^1\in\sigma_1}\frac{|v_n|^2}{z_n^1-\mu_p}=\alpha G'(\mu_p)(\not=0)$$
where $G(z)$ is given by \eqref{eq3.12} and $G'(\mu_p)>0$. Define the vectors
\begin{equation}
\widetilde{u}_p\stackrel{\rm def}{=}\frac1{\sqrt{G'(\mu_p)}}\sum\limits_{z_n^1\in\sigma_1}\frac{v_n}{z_n^1-\mu_p}u_n.\label{eq3.15}
\end{equation}
This series converges in $H$ and $\|\widetilde{u}_p\|=1$. Moreover, $\widetilde{u}_p$ are eigenfunctions of the operator $L_{\alpha,1}$, $L_{\alpha,1}\widetilde{u}_p=\mu_p\widetilde{u}_p$. Really,
$$L_{\alpha,1}\widetilde{u}_p=\frac1{\sqrt{G'(\mu_p)}}\sum\limits_{z_n^1\in\sigma_1}\frac{v_n}{z_n^1-\mu_p}z_n^1u_n+\frac\alpha{\sqrt{G'(\mu_p)}}\sum\limits_{z_n^1\in\sigma_1}\frac{|v_n|^2}{z_n^1-\mu_p}v$$
$$=\mu_p\widetilde{u}_p+\frac1{\sqrt{G'(\mu_p)}}Q(\mu_p)v=\mu_p\widetilde{u}_p$$
due to $Q(\mu_p)=0$. Formula \eqref{eq3.15} implies that $c_p(f)$ \eqref{eq3.14} equal $c_p(f)=\langle f,\widetilde{u}_p\rangle\widetilde{u}_p$. Therefore resolvent $R_{L_{\alpha,1}}(z)$ \eqref{eq3.11} is given by
\begin{equation}
R_{L_{\alpha,1}}(z)f=\sum\limits_{\mu_p\in\sigma(L_{\alpha,1})}\frac{\widetilde{f}_p}{\mu_p-z}\widetilde{u}_p\label{eq3.16}
\end{equation}
where $\widetilde{f}_p=\langle f,\widetilde{u}_p\rangle$ are Fourier coefficients of the vector $f$ in the basis $\{\widetilde{u}_p\}$ \eqref{eq3.15} of the subspace $H_1$.

\begin{lemma}\label{l3.4}
Spectrum of the operator $L_{\alpha,1}$ \eqref{eq3.9} is simple and coincides with the set $\sigma(L_{\alpha,1})$ \eqref{eq3.13}. Eigenvectors of $L_{\alpha,1}$ corresponding to eigenvalues $\mu_p$ are given by \eqref{eq3.15} and resolvent $R_{L_{\alpha,1}}(z)$ is given by formula \eqref{eq3.16}.
\end{lemma}

It can turn out that a number $\mu_p\in\sigma(L_{\alpha,1})$ coincides with $z_n^0\in\sigma_0$ \eqref{eq3.7}, then the proper subspace of the operator $L_\alpha$ \eqref{eq3.3} is two-dimensional and is generated by the vectors $\{u_n,\widetilde u_p\}$.

\begin{theorem}\label{t3.1}
Spectrum of the operator $L_\alpha$ \eqref{eq3.3} is
\begin{equation}
\sigma(L_\alpha)=\sigma_0\cup\sigma(L_{\alpha,1})\label{eq3.17}
\end{equation}
where $\sigma_0$ and $\sigma(L_{\alpha,1})$ are given by \eqref{eq3.7} and \eqref{eq3.13} correspondingly. Spectrum of the operator $L_\alpha$ is simple, excluding points belonging to $\sigma_0\cap\sigma(L_{\alpha,1})$ where it is of multiplicity 2.
\end{theorem}

\begin{remark}\label{r3.1}
Let the set $\sigma(L_0)$ \eqref{eq3.1} be such that
\begin{equation}
\lim\limits_{n\rightarrow\infty}z_n=\infty;\quad\min\limits_{n,s}|z_n-z_s|=d>0,\label{eq3.18}
\end{equation}
then the intersection $\sigma_0\cap\sigma(L_{\alpha,1})$ is finite. Thus, when conditions \eqref{eq3.18} hold, the operator $L_\alpha$ \eqref{eq3.3} has finite number of points of the spectrum of multiplicity $2$.
\end{remark}

Really, if $\sigma_0\cap\sigma(L_{\alpha,1})$ is infinite and consists of points $\{\widehat{z}_p^0\}_1^\infty$, then in view of $Q(\widehat{z}_p^0)=0$
$$1+\alpha\sum\limits_{z_n^1\in\sigma_1}\frac{|v_n|^2}{z_n^1-\widehat{z}_p^0}=0\quad(\forall p).$$
Formula \eqref{eq3.18} and $\sum|v_n|^2=1$ imply that this series converges uniformly. Passing to the limit as $p\rightarrow\infty$ in this equality and taking into account \eqref{eq3.18}, one obtains $1=0$ which is absurd.

Lemma \ref{l1.4} implies that the conditions \eqref{eq3.18} for the operator $L_0$ \eqref{eq1.1}, \eqref{eq1.2} hold, therefore for the operator $L_\alpha$ \eqref{eq2.1} Theorem \ref{t3.1} holds and there is finite number of points in $\sigma(L_\alpha)$ of multiplicity 2.
\vspace{5mm}

{\bf 3.3.} Calculate $Q(z)=1+\alpha\langle R_{L_0}(z)v,v\rangle$ for the operator $L_0$ \eqref{eq1.1}, \eqref{eq1.2}. Relation \eqref{eq1.32} implies that
$$\langle R_{L_0}(\lambda^3)v,v\rangle=\frac i{\lambda^4\Delta(0,\lambda)}\left\{\langle\int\limits_0^l[d(i\lambda(x-l))d(-i\lambda t)-d(i\lambda x)d(i\lambda(l-t))]\right.$$
$$\times v(t)dt,v(x)\rangle+\left\langle d(i\lambda l)\int\limits_0^xd(i\lambda(x-t))v(t)dt,v(x)\right\rangle$$
$$\left.-\left\langle d(-i\lambda l)\int\limits_x^ld(i\lambda(x-t))v(t)dt,v(x)\right\rangle\right\}.$$
Taking into account \eqref{eq2.2}, \eqref{eq2.5} and
$$m(\lambda)=\left\langle\int\limits_0^xd(i\lambda(x-t)v(t)dt,v(x)\right\rangle,$$
one obtains
$$\langle R_{L_0}(\lambda^3)v,v\rangle=\frac i{\lambda^4\Delta(0,\lambda)}\{m(\lambda)d(i\lambda l)-\widetilde{v}_d^*(\lambda)\widetilde{w}_d(-\lambda)$$
$$-m^*(\lambda)d(-i\lambda l)+\widetilde{v}_d(\lambda)\widetilde{w}_d^*(-\lambda)\}=\frac i{\lambda^4\Delta(0,\lambda)}\{F(\lambda)-F^*(\lambda)\}$$
where $F(\lambda)$ is given by \eqref{eq2.15}, using \eqref{eq2.14}, one arrives at the equality
$$\alpha\langle R_{L_0}(\lambda^3)v,v\rangle=\frac{\Delta(\alpha,\lambda)-\Delta(0,\lambda)}{\Delta(0,\lambda)}.$$

\begin{lemma}\label{l3.5} For the operator $L_0$ \eqref{eq1.1}, \eqref{eq1.2}, the identity
\begin{equation}
Q(\lambda^3)=\frac{\Delta(\alpha,\lambda)}{\Delta(0,\lambda)}\label{eq3.19}
\end{equation}
where $Q(z)=1+\alpha\langle R_{L_0}(z)v,v\rangle$ and $\Delta(0,\lambda)$ \eqref{eq1.20} and $\Delta(\alpha,\lambda)$ \eqref{eq2.14} are characteristic functions of the operators $L_0$ \eqref{eq1.1}, \eqref{eq1.2} and $L_\alpha$ \eqref{eq2.1}.
\end{lemma}

\section{Inverse problem}\label{s4}

{\bf 4.1.} To obtain multiplicative expansions of characteristic functions $\Delta(0,\lambda)$ \eqref{eq1.20} and $\Delta(\alpha,\lambda)$ \eqref{eq2.14}, use the well-known Hadamard theorem on factorization \cite{9, 10}.

\begin{theorem}[Hadamard]\label{t4.1}
Let $f(z)$ be an entire function of exponential type and $\{z_n\}_1^\infty$ be a subsequence of its zeros where each zero is repeated according to its multiplicity, besides, $z_1\not=0$ and $0<|z_1|\leq|z_2|\leq...$. Then
\begin{equation}
f(z)=cz^ke^{bz}\prod\limits_{n}\left(1-\frac z{z_n}\right)e^{z/z_n}\label{eq4.1}
\end{equation}
where $k\in\mathbb{Z}_+$; $c$, $b\in\mathbb{C}$.
\end{theorem}

Lemma \ref{l1.4} implies that zeros of the function $\Delta(0,\lambda)$ come in triplets $\{\zeta_l\lambda_n(0)\}_{l=1}^3$ (and $\{-\zeta_l\lambda_n(0)\}_{l=1}^3$) ($\lambda_n(0)>0$), therefore
\begin{equation}
\prod\limits_{n=1}^3\left(1-\frac\lambda{\zeta_k\lambda_n(0)}\right)e^{\lambda/\zeta_k\lambda_n(0)}=1-\frac{\lambda^3}{\lambda_n^3(0)}\label{eq4.2}
\end{equation}
due to \eqref{eq1.4} (similarly for $\{-\zeta_l\lambda_n(0)\}$ also), and according to \eqref{eq4.1}
\begin{equation}
\Delta(0,\lambda)=c\lambda^ke^{b\lambda}\prod\limits_n\left(1-\frac{\lambda^6}{\lambda_n^6(0)}\right).\label{eq4.3}
\end{equation}
Infinite product \eqref{eq4.3} converges uniformly on every compact due to location of $\lambda_n(0)$ (Lemma \ref{l1.4}). Taylor formula (ix) \eqref{eq1.10} for $d(z)$ and \eqref{eq1.20} imply the expansion
\begin{equation}
\Delta(0,\lambda)=-l^2+\frac28\lambda^5l^8+...\label{eq4.4}
\end{equation}
and thus $\Delta(0,0)=-l^2$ ($\not=0$) and $k=0$, $c=-l^2$ in the formula \eqref{eq4.3}. So,
$$\Delta(0,\lambda)=-l^2e^{b\lambda}\prod\limits_n\left(1-\frac{\lambda^6}{\lambda_n^6(0)}\right).$$
Differentiating this equality and assuming that $\lambda=0$, one obtains that $b=0$ since $\Delta'(0,0)=0$ due to \eqref{eq4.4}. Thus, multiplicative expansion of $\Delta(0,\lambda)$ is
\begin{equation}
\Delta(0,\lambda)=-l^2\prod\limits_n\left(1-\frac{\lambda^6}{\lambda_n^2(0)}\right).\label{eq4.5}
\end{equation}

Property \eqref{eq2.16} of the characteristic function $\Delta(\alpha,\lambda)$ \eqref{eq2.14} yields that its roots also form triplets $\{\zeta_l\lambda_n(\alpha)\}_{l=1}^3$ ($0=\lambda_n(\alpha)\in\mathbb{R}$) and for them equality \eqref{eq4.2} is true, therefore, due to Theorem \ref{t4.1},
\begin{equation}
\Delta(\alpha,\lambda)=c\lambda^ke^{b\lambda}\prod\limits_n\left(1-\frac{\lambda^3}{\lambda_n^3(\alpha)}\right).\label{eq4.6}
\end{equation}
This product also converges uniformly on every compact. Again using the Taylor formula (ix) \eqref{eq1.10} for $d(z)$, one finds series expansion of the function $F(\lambda)$ \eqref{eq2.15},
\begin{equation}
F(\lambda)=\frac{\lambda^4}4a_1(v)-\frac{i\lambda^7}{2\cdot5!}a_2(v)+...\label{eq4.7}
\end{equation}
where
\begin{equation}
a_1(v)=\left\langle l^2\int\limits_0^x(x-t)^2v(t)dt+(l-x)^2\int\limits_0^lt^2v(t)dt,v(x)\right\rangle.\label{eq4.8}
\end{equation}
Substituting \eqref{eq4.4} and \eqref{eq4.7} into \eqref{eq2.14}, one finds the expansion
\begin{equation}
\Delta(\alpha,\lambda)=-l^2+{\alpha i}4(a_1(v)-\overline{a_1(v)})+\frac{\lambda^3\alpha}{2\cdot 5!}(a_2(v)-\overline{a_2(v)})+...\label{eq4.9}
\end{equation}
Hence \eqref{eq4.6} implies that $k=0$ and
\begin{equation}
c=-l^2+\frac{\alpha i}4(a_1(v)-\overline{a_1(v)})\not=0.\label{eq4.10}
\end{equation}
Take into account that $\Delta'(\alpha,0)=0$ \eqref{eq4.9}, then \eqref{eq4.6} implies that $b=0$. So,
\begin{equation}
\Delta(\alpha,\lambda)=c\prod\limits_n\left(1-\frac{\lambda^3}{\lambda_n^3(\alpha)}\right)\label{eq4.14}
\end{equation}
where $c$ is given by \eqref{eq4.10}. Further, the number $c$ will be found in terms of spectral data.
\vspace{5mm}

{\bf 4.2.} Use equality \eqref{eq3.19},
\begin{equation}
Q(z)=\frac{\Delta(\alpha,z^{1/3})}{\Delta(0,z^{1/3})},\label{eq4.12}
\end{equation}
and note that in the quotient all terms in infinite products \eqref{eq4.5}, \eqref{eq4.14} which correspond to $z_n(0)=\pm\lambda_n^3(0)$ from $\sigma_0$ \eqref{eq3.7} reduce. Only terms corresponding to $z_n(0)=\varepsilon\lambda_n^3(0)\in\sigma_1$ \eqref{eq3.7} ($\varepsilon=1$ or $-1$) and corresponding to it in numerator $z_n(\alpha)$ remain. So,
\begin{equation}
1+\alpha\sum\limits_{z_n\in\sigma_1}\frac{|v_p|^2}{z_n-z}=-\frac c{l^2}\prod\limits_{z_n\in\sigma_1}\frac{z_n}{z_n(\alpha)}\left(1+\frac{z_n(\alpha)-z_n}{z_n-z}\right).\label{eq4.13}
\end{equation}
Substituting $z=iy$ ($y\in\mathbb{R}$) and passing to the limit as $y\rightarrow\infty$, one obtains
$$1=-\frac c{l^2}\prod\limits_{z_n\in\sigma_1}\frac{z_n}{z_n(\alpha)},$$
hence once finds the number $c$ via the spectral data
\begin{equation}
c=-l^2\prod\limits_{z_n\in\sigma_1}\left(1+\frac{z_n(\alpha)-z_n}{z_n}\right).\label{eq4.14}
\end{equation}
Calculating residue at the point $z=z_n$ in both parts of the equality \eqref{eq4.13}, one has
\begin{equation}
\alpha|v_p|^2=-\frac c{l^2}\frac{z_p}{z_p(\alpha)}(z_p(\alpha)-z_p)\prod\limits_{n\not=p}\frac{z_n}{z_n(\alpha)}\left(1+\frac{z_n(\alpha)-z_n}{z_n-z_p}\right).\label{eq4.15}
\end{equation}

\begin{theorem}\label{t4.2}
Using the formulas \eqref{eq4.15}, the numbers $\alpha|v_p|^2$, where $c$ is given by \eqref{eq4.14}, are found unambiguously from spectra $\sigma(L_0)=\{z_n\}$ and $\sigma(L_\alpha)=\{z_n(\alpha)\}$ of the operators $L_0$ \eqref{eq1.1}, \eqref{eq1.2} and $L_\alpha$ \eqref{eq2.1}.
\end{theorem}

\begin{remark}\label{r4.1}
The number $\alpha$ is defined by the numbers $\alpha|v_p|^2$ \eqref{eq4.15} from the condition $\sum|v_p^2|=1$. Thus, $\{\alpha,|v_p|^2\}$ is unambiguously recovered by the spectra $\sigma(L_0)$ and $\sigma(L_\alpha)$.
\end{remark}

Relation \eqref{eq1.25} implies that Fourier coefficient equals
\begin{equation}
\begin{array}{ccc}
{\displaystyle v_p=\langle v(x),u(0,\lambda_p(0),x)\rangle=\frac1{u_p}\{s(-i\lambda_p(0)l)\widetilde{v}_s(\lambda_p(0))}\\
+(1-c(-i\lambda_p(0)l))\widetilde{v}_d(\lambda_p(0))\}
\end{array}\label{eq4.16}
\end{equation}
where $\widetilde{v}_s(\lambda)$ and $\widetilde{v}_d(\lambda)$ are given by \eqref{eq2.2}. For the function $g(x)=l-x$, Fourier coefficients are calculated explicitly,
\begin{equation}
g_p=\langle g(x), u(0,\lambda_p(0),x)\rangle=\frac{il}{\lambda_p(0)}(c(-i\lambda_p(0))-1).\label{eq4.17}
\end{equation}

\begin{theorem}\label{t4.3}
The number $\alpha$ and the function $v(x)$ ($\|v\|=1$) are unambiguously recovered from the spectra $\sigma(L_0)$, $\sigma(L_\alpha(v))$, $\sigma(L_\alpha(v+g))$, $\sigma(L_\alpha(v+ig))$ where $L_\alpha(v)$ is given by \eqref{eq2.1} and $g(x)=l-x$.
\end{theorem}

P r o o f. Theorem \ref{t4.2} implies that numbers $\alpha|v_p|^2$ are unambiguously defined by $\sigma(L_0)$ and $\sigma(L_\alpha(v))$. Exactly in the same way from $\sigma(L_0)$ and $\sigma(L_\alpha(v+g))$ one finds
$$\alpha|v_p+g_p|^2=\alpha|v_p|^2+2\Re(\alpha v_p\overline{g_p})+|g_p|^2$$
where $g_p$ are given by \eqref{eq4.17}. Hence, the numbers $\Re(\alpha v_p\overline{g_p})$ are unambiguously calculated by the three spectra $\sigma(L_0)$, $\sigma(L_\alpha(v))$, $\sigma(L_\alpha(v+g))$. Similarly, from $\sigma(L_0)$, $\sigma(L_\alpha(v))$, $\sigma(L_\alpha(v+ig))$ $\Im(\alpha v_p\overline{g_p})$ are defined. Thus, $\alpha v_p\overline{g_p}$, and so $\alpha v_p$ also, are unambiguously calculated from the four spectra $\sigma(L_0)$, $\sigma(L_\alpha(v))$, $\sigma(L_\alpha(v+g))$, $\sigma(L_\alpha(v+ig))$. Finally, numbers $\alpha$ and $v_p$ are found from the normalization condition $\sum|v_p|^2=1$. Thereafter, the function $v(x)$ is defined by its Fourier series,
$$v(x)=\sum\limits_pv_pu(0,\lambda_p(0),x).\quad\blacksquare$$
\vspace{5mm}

{\bf 4.3.} Proceed to description of the data of inverse problem. Denote by $\sigma_1$ a subset in $\sigma(L_0)=\{z_n(0,\varepsilon)=\varepsilon\lambda_n^3(0):\Delta(0,\lambda_n(0))=0,\varepsilon=\pm1\},$
\begin{equation}
\sigma_1\stackrel{\rm def}{=}\{z_n^1(0,\varepsilon):z_n^1(0,\varepsilon)\in\sigma(L_0)\}\label{eq4.18}
\end{equation}
where the numbers $z_n^1(0,\varepsilon)$ are enumerated in ascending order (the case of $\sigma_1=\sigma(L_0)$ is not excluded). And let
\begin{equation}
\sigma_0\stackrel{\rm def}{=}\sigma(L_0)\backslash\sigma_1=\{z_n^0(0,\varepsilon):z_n^0(0,\varepsilon)=\sigma(L_0)\}.\label{eq4.19}
\end{equation}
Consider the set
\begin{equation}
\sigma(\mu)\stackrel{\rm def}{=}\{\mu_n\in\mathbb{R}:n\in\mathbb{N}\}\label{eq4.20}
\end{equation}
assuming that sets $\sigma(\mu)$ and $\sigma_1$ interlace, moreover, $\mu_k$ are enumerated in ascending order, $\mu_k\not=\mu_s$ ($k\not=s$), $\mu_1\not=0$ ($\sigma(\mu)\cap\sigma(1)=\emptyset$).

\begin{condition}\label{c1}
The intersection $\sigma(\mu)\cap\sigma_0$ is finite.
\end{condition}

Since the number $c$ \eqref{eq4.14} is finite, then one arrives at the second condition.

\begin{condition}\label{c2}
The sets $\sigma_1$ and $\sigma(\mu)$ are such that the series
\begin{equation}
\sum\limits_n\ln\left(1+\frac{\mu_n-z_n^1(0,\varepsilon)}{z_n^1(0,\varepsilon)}\right)<\infty\label{eq4.21}
\end{equation}
which is equivalent to the convergence of the product
\begin{equation}
\prod\limits_n\left(1+\frac{\mu_n-z_n^1(0)}{z_n^1}\right).\label{eq4.22}
\end{equation}
\end{condition}

Numbers $v_p\in l^2({\mathbb N})$, therefore the following requirement is natural, due to \eqref{eq4.15}.

\begin{condition}\label{c3}
The numerical series
\begin{equation}
\sum\limits_p\frac{z_p^1(0,\varepsilon)}{\mu_p}(\mu_p-z_p^1(0,\varepsilon))\left[\prod\limits_{n\not=p}\frac{z_n^1(0,\varepsilon)}{\mu_n}\left(1+\frac{\mu_n-z_n^1(0,\varepsilon)}{z_n^1(0)-z_p^1(0,\varepsilon)}\right)\right]<\infty
\label{eq4.23}
\end{equation}
converges.
\end{condition}

Form the infinite product
\begin{equation}
b_{\sigma_1}(z)\stackrel{\rm def}{=}\prod\limits_{z_n^1\in\sigma_1}\left(1-\frac z{z_n^1(0,\varepsilon)}\right)\label{eq4.24}
\end{equation}
converging in any circle $C_R=\{z\in\mathbb{C}:|z|<R\}$ due to $z_n\in\sigma(L_0)$. Similarly, set
\begin{equation}
b_{\sigma(\mu)}(z)\stackrel{\rm def}{=}\prod\limits_{z_n^1\in\sigma_1}\left(1-\frac z{\mu_n}\right).\label{eq4.25}
\end{equation}
Consider a meromorphic function
$$Q(z)=k\frac{b_{\sigma(\mu)}(z)}{b_{\sigma_1}(z)}\quad(k\in\mathbb{R}).$$
Residue of $Q(z)$ at the point $z_p^1(0,\varepsilon)$ equals
\begin{equation}
c_p=\lim\limits_{z\rightarrow z_p^1(0)}(z_p^1(0,\varepsilon)-z)Q(z)=kz_p^1(0)\frac{b_{\sigma(\mu)}(z_n^1(0,\varepsilon))}{\displaystyle\prod\limits_{n\not=p}\left(1-\frac{z_p^1(0,\varepsilon)}{z_n^1(0,\varepsilon)}\right)}.\label{eq4.26}
\end{equation}

\begin{remark}\label{r4.2}
Interlacing of sequences $z_n^1(0)$ and $\mu_n$ implies that all the numbers $c_p$ are of the same sign, $\sign c_p=\sign c_1$ ($\forall p$).
\end{remark}

If Condition \ref{c3} is met, then defining the number
$$\frac1\alpha\stackrel{\rm def}{=}\sum\limits_p\frac{z_p^1(0,\varepsilon)}{\mu_p}(\mu_p-z_p^1(0,\varepsilon))\prod\limits_{n\not=p}\frac{z_n^1(0,\varepsilon)}{\mu_n}\left(1+\frac{\mu_n-z_n^1(0,\varepsilon)}{z_n^1(0)-z_p^1(0,\varepsilon)}\right)$$
one obtains that the numbers
$$|v_p|^2\stackrel{\rm def}{=}\frac{c_p}\alpha\quad(p\in\mathbb{N})$$
are positive and $\sum|v_p|^2=1$.

\begin{theorem}\label{t4.4}
Let $\sigma$ be a countable set on $\mathbb{R}$ and $\sigma=\sigma_0\cup\sigma(\mu)$ where $\sigma_0\subset\sigma(L_0)$ \eqref{eq4.19} ($\sigma(L_0)=\{z_n(0,\varepsilon)=\varepsilon\lambda_n^3(0):\Delta(0,\lambda_n(0))=0,\varepsilon=\pm1\}$ and $\sigma(\mu)$ \eqref{eq4.20} consists of pairwise different points $\mu_n$ ($\mu_1\not=0$), $n\in\mathbb{N}$, enumerated in ascending order and interlacing with $z_n^1(0,\varepsilon)\in\sigma_1=\sigma(L_0)\backslash\sigma_0$ and $\sigma_1\cap\sigma(\mu)=\emptyset$.

In order that $\sigma$ coincide with the spectrum $\sigma(L_\alpha)$ of the operator $L_\alpha$ \eqref{eq2.1} it is necessary and sufficient that

(a) $\sigma(\mu)\cap\sigma_0$ is finite;

(b) the series
\begin{equation}
\begin{array}{ccc}
{\displaystyle\sum\ln\left(1+\frac{\mu_n-z_n^1(0,\varepsilon)}{z_n^1(0,\varepsilon)}<\infty\right)};\\
{\displaystyle\sum\limits_p\frac{z_p^1(0,\varepsilon)}{\mu_p}(\mu_p-z_n^1(0,\varepsilon))\prod\limits_{n\not=p}\left(\frac{\mu_n-z_p^1(0,\varepsilon)}{z_n^1(0)-z_p^1(0,\varepsilon)}\right)<\infty}
\end{array}\label{eq4.26}
\end{equation}
converge.
\end{theorem}

Description of the spectrum $\sigma(L_\alpha)$ of the operator $L_\alpha$ \eqref{eq2.1} one can give in terms of characterization of the class of functions $\Delta(\alpha,\lambda)$ that are characteristic for $L_\alpha$ \cite{7}.

\renewcommand{\refname}{ \rm

\end{Large}
\end{document}